\documentclass[12pt]{article}
\usepackage{color}
\definecolor{darkblue}{rgb}{0.00,0.25,0.50}
\usepackage[colorlinks,filecolor=blue,citecolor=darkblue]{hyperref}

\usepackage{amssymb,amsmath}
\usepackage{url}
\usepackage{amscd, color}
\usepackage[english]{babel}
\begin{document}

\thispagestyle{empty}

\begin{center}
\textbf{UNIFORM APPROXIMATION OF POISSON INTEGRALS OF FUNCTIONS FROM THE CLASS $
 H_\omega $  BY DE LA VALL\'{E}E POUSSIN
SUMS}
\end{center}
\vskip0.5cm
\begin{center}
A.S. SERDYUK and Ie.Yu. OVSII
\end{center}
\vskip0.5cm

\abstract{We obtain asymptotic equalities for least upper bounds of deviations in the
uniform metric of de la Vall\'{e}e Poussin sums on the sets $C^{q}_{\beta}H_\omega$ of
Poisson integrals of functions from the class $H_\omega $ generated by convex upwards moduli
of continuity $\omega (t)$ which satisfy the condition $\omega (t)/t\to\infty$ as $t\to0.$
As an implication, a solution of the Kolmogorov--Nikol'skii problem for de la Vall\'{e}e
Poussin sums on the sets of Poisson integrals of functions belonging to Lipschitz classes
$H^\alpha ,$ $0<\alpha <1,$ is obtained.}



\section{Introduction}

Let $L$ be the space of $2\pi $-periodic summable functions $f(t)$ with the norm
\mbox{$\|f\|_{L}=\int\limits_{-\pi }^{\pi}{
    |f(t)|}\,dt,$}
let $L_\infty$ be the space of $2\pi $-periodic functions $f(t)$ with the norm
\mbox{$\|f\|_{L_\infty}=\mathop{\rm ess\,sup}\limits_{t\ }|f(t)|,$} and let $C$ be the space
of continuous $2\pi $-periodic functions $f(t)$ in which the norm is defined by the formula
$\|f\|_C=\max\limits_{t}|f(t)|.$

Further, let $C^q_\beta \mathfrak{N} $ be the set of the Poisson integrals of functions
$\varphi $ from $\mathfrak{N}\subset L,$ i.e., functions $f$ of the form
\begin{equation}\label{6.05.10-14:46:39}
    f(x)=A_0+\frac{1}{\pi }\int_{0}^{2\pi }\varphi (
    x+t)P_{q,\beta }(t)\,dt,\ \ A_0\in \mathbb{R},\ \ x\in \mathbb{R},\ \ \varphi \in
    \mathfrak{N},
\end{equation}
where
    $$P_{q,\beta }(t)=
    \sum\limits_{k=1}^{\infty}q^k\cos \left(kt+\frac{\beta \pi }{2}\right),\ \
    q\in (0,1),\ \ \beta \in \mathbb{R},$$
is the Poisson kernel with parameters $q$ and $\beta .$ If $\mathfrak{N}=U_\infty$, where
$$U_\infty=\{\varphi \in L_\infty: \|\varphi \|_{L_\infty}\leqslant 1\},$$ then the classes
$C^q_\beta \mathfrak{N} $ will be denoted by $C^q_{\beta ,\infty},$ and if
$\mathfrak{N}=H_\omega $, where
    $$H_\omega =\{\varphi \in C:\ |\varphi (t')-\varphi (t'')|
    \leqslant \omega (|t'-t''|)\ \ \forall t',t''\in \mathbb{R}\},$$
and $\omega (t)$ is an arbitrary modulus of continuity, then $C^q_\beta \mathfrak{N} $ will
be denoted by $C^q_\beta H_\omega . $

Denoting by $S_k(f;x)$ the $k$th partial sum of the Fourier series of the summable function
$f$, we associate each function $f\in C^q_\beta H_\omega $ with the trigonometric polynomial
of the form
\begin{equation}\label{6.05.10-14:55:11}
    V_{n,p}(f;x)=\frac{1}{p}\sum\limits_{k=n-p}^{n-1}S_k(f;x),\ \ p,n\in
    \mathbb{N},\ \ p\leqslant n.
\end{equation}
The sums $V_{n,p}(f;x)$ appeared in 
\cite{Vallee
Poussin} and are called the de la Vall\'{e}e Poussin sums with parameters $n$ and $p$.

The purpose of the present work is to solve the Kolmogorov--Nikol'skii problem for de la
Vall\'{e}e Poussin sums, which consists of obtaining an asymptotic equality as
\mbox{$n-p\to\infty$} of the quantity
\begin{equation}\label{6.05.10-15:01:20}
    \mathcal{E}(\mathfrak{N};V_{n,p})=\sup\limits_{f\in \mathfrak{N}
    }{\|f(\cdot)-V_{n,p}(f;\cdot)\|_C,}
\end{equation}
where $\mathfrak{N}=C^q_\beta H_\omega ,$ $q\in(0,1),$ $\beta \in \mathbb{R}$ and $\omega
(t)$ is a given convex upwards modulus of continuity. Since for $p=1$
$V_{n,p}(f;x)=V_{n,1}(f;x)=S_{n-1}(f;x),$ the quantity
    $$\mathcal{E}(\mathfrak{N};S_{n-1})=\sup\limits_{f\in\mathfrak{N}}{\|f(\cdot)-
    S_{n-1}(f;\cdot)\|_C}$$
is the special case of (\ref{6.05.10-15:01:20}).

The problem of obtaining asymptotic equalities for the quantities of the form
    $$\sup\limits_{f\in\mathfrak{N}}{\|f(\cdot)-V_{n,p}(f;\cdot)\|_X},$$
has a rich history in various function classes $\mathfrak{N}$ and metrics $X\subset L$ and
connected with the names of {\sc A.N. Kolmogorov} \cite{Kolmog}, {\sc S.M. Nikol'skii}
\cite{NIKOLSKIY-1946-1}, {\sc A.F. Timan} \cite{TIMAN-60}, {\sc \mbox{S.B. Stechkin}}
\cite{Stechkin_1980}, {\sc A.V. Efimov} \cite{EFIMOV60}, {\sc S.A. Telyakovskii}
\cite{TELYKOVSKIY-61}, {\sc A.I. Stepanets} \cite{Stepanets_2001}, {\sc V.I. Rukasov}
\cite{RUKASOV-2003} and many others. See  \cite{Rukasov_Chaichenko}, \cite{Serdyuk V-P},
\cite{Serdyuk V-P-2009}, \cite{Serdyuk V-P-2010}, \cite{STEANETS-RUKASOV} and
\cite{STEANETS-RUKASOV-CHAICHENKO} for more details on the history of this problem.

{\sc S.M. Nikol'skii} \cite{NIKOLSKIY-1946-1} proved that, if $n\to\infty$ then
\begin{equation}\label{9.11.10-14:52:26}
    \mathcal{E}(C^q_{\beta ,\infty};S_{n-1})=q^n
    \bigg(\frac{8}{\pi ^2}\textbf{K}(q)+O(1)n^{-1}\bigg), \ \ \beta \in \mathbb{R},
\end{equation}
where
    $$\textbf{K}(q):=\int_{
    0}^{\pi /2}\frac{dt}{\sqrt{1-q^2\sin^2 t}},\ \ q\in(0,1)$$
is the complete elliptic integral of the 1st kind. In \cite{Stechkin_1980} {\sc S.B.
Stechkin} improved the estimate of the remainder term in (\ref{9.11.10-14:52:26}) by showing
that
    \begin{equation}\label{15.02.11-11:19:44}
    \mathcal{E}(C^q_{\beta ,\infty};S_{n-1})=q^n
    \left(\frac{8}{\pi ^2}\textbf{K}(q)+O(1)\frac{q}{(1-q)n}\right), \ \ \beta \in \mathbb{R},
    \end{equation}
where $O(1)$ is a quantity uniformly bounded in $n,$ $q$ and $\beta .$

In studying approximate properties of the Fourier sums $S_{n-1}(f;x)$ on the classes
$C^q_\beta H_\omega $, {\sc A.I. Stepanets} \cite{Stepanets_2001} proved that for any
$q\in(0,1),$ $\beta \in\mathbb{R}$ and arbitrary modulus of continuity $\omega (t)$ the
equality
\begin{equation}\label{11.05.10-12:13:35}
    \mathcal{E}(C^q_\beta H_\omega ;S_{n-1})=q^n\bigg(\frac{4}{
    \pi^2 }\textbf{K}(q)e_n(\omega )+O(1)\frac{\omega (1/n)}{(1-q)^2n}\bigg),
\end{equation}
holds as $n\to\infty$, where
\begin{equation}\label{22.12.10-23:18:36}
    e_n(\omega ):=\theta_\omega\int_{
    0}^{\pi /2}\omega \left(\frac{2t}{n}\right)\sin t\,dt,
\end{equation}
$\theta_{\omega} \in [1/2,1]$,  ($\theta_\omega =1$ if $\omega (t)$ is a convex upwards
modulus of continuity) and $O(1)$ is the same as in (\ref{15.02.11-11:19:44}).

Reasoning as in \cite{Stepanets_2001}, {\sc V.I. Rukasov} and {\sc S.O. Chaichenko}
\cite{Rukasov_Chaichenko} showed that if \mbox{$q\in(0,1),$} $\beta \in\mathbb{R},$ $n,p\in
\mathbb{N},$ $p\leqslant n$ and $\omega (t)$  is an arbitrary modulus of continuity, then as
$n\to\infty,$
\begin{equation}\label{7.05.10-12:08:22}{
    \mathcal{E}(C^q_\beta H_\omega ;V_{n,p})=
    \frac{2q^{n-p+1}}{\pi (1-q^2)p}e_{n-p+1}(\omega )+}$$
    $${+O(1)\frac{q^{n-p+1}}{p}\omega \bigg(\frac{1}{n-p+1}\bigg)\bigg(\frac{q^{p}}{
    1-q^2}+\frac{1}{(1-q)^3(n-p+1)}\bigg),}
\end{equation}
where $O(1)$ is a quantity uniformly bounded in $n,$ $p,$ $q$ and $\beta .$

For the classes $C^q_\beta H_\omega $ there are well-known estimates of the best uniform
approximations by trigonometric polynomials of order not more than $\leqslant n-1$ (see, for
example, \cite[p. 509]{Stepanets_2002}):
    $$E_n(C^q_{\beta}H_\omega )=\sup\limits_{f\in C^q_{\beta}H_\omega }{
    \inf\limits_{t_{n-1}}{\|f(\cdot)-t_{n-1}(\cdot)\|_C}}\asymp q^n\omega (1/n)$$
(the notation $\alpha (n)\asymp\beta (n)$ as $n\to\infty$ means that there exist $K_1,
K_2>0$ such that $K_1\alpha (n)\leqslant \beta (n)\leqslant K_2\alpha (n)$).  It's easy to
see from (\ref{7.05.10-12:08:22}) that if $n\to\infty$ and the values of the parameter $p$
are bounded, then the de la Vall\'{e}e Poussin sums realize the order of the best uniform
approximation on the classes $C^q_\beta H_\omega $. But in that case estimate
(\ref{7.05.10-12:08:22}) isn't an asymptotic equality, taking the form
    $$\mathcal{E}(C^q_\beta H_\omega ;V_{n,p})=
    O(1)\frac{q^{n-p+1}}{p(1-q^2)}\omega \bigg(\frac{1}{n-p+1}\bigg).$$

Thus for $n-p\to\infty$ and bounded values of the parameter $p$ the question of the
asymptotic behavior of a principal term of the quantities $\mathcal{E}(C^q_\beta H_\omega
;V_{n,p})$ was still unknown.

The proofs of (\ref{11.05.10-12:13:35}) and (\ref{7.05.10-12:08:22}) are based on the
well-known Korneichuk--Stechkin lemma (see, e.g., \cite{Stepanets_2001}). In the present
work we use a somewhat different way which, as will be seen later, proves itself in some
important cases.

\section{Main result}
\hypertarget{15.04.11-16:12:48}{}

The main result of the paper is the next theorem.

\textbf{Theorem {1}.}  \emph{Let $q\in(0,1),$ $\beta \in\mathbb{R},$ $n,p\in \mathbb{N},$
$p<n$ and let $\omega (t)$ be a convex upwards modulus of continuity. Then as
$n-p\to\infty$$:$
\begin{equation}\label{27.03.10-15:23:27}
    \mathcal{E}(C^q_\beta H_\omega ;V_{n,p})=
    \frac{q^{n-p+1}}{p}\Bigg(
    \frac{{{K}}_{p,q}}{\pi ^2}e_{n-p+1}(\omega )
    +O(1)\frac{\omega (\pi )}{(1-q)^{\delta (p)}(n-p+1)}\Bigg),
\end{equation}
where
\begin{equation}\label{29.03.10-12:30:45}
    {{K}}_{p,q}:=\int_{0}^{2\pi }\frac{\sqrt{1-2q^p\cos pt+q^{2p}}}{
    1-2q\cos t+q^2}\,dt,
    \end{equation}
    $$e_{n-p+1}(\omega )=\int_{
    0}^{\pi /2}\omega \left(\frac{2t}{n-p+1}\right)\sin t\,dt,$$
    \begin{equation}\label{29.03.10-12:31:28}
        \delta (p):=
  \begin{cases}
    2, & p=1, \\
    3, & p=2,3,\ldots,
  \end{cases}
    \end{equation}
and $O(1)$ is a quantity uniformly bounded in $n,$ $p,$ $q$, $\omega $ and $\beta$.}

Since
    $$\frac{2}{\pi }\omega (\frac{\pi }{k})\leqslant e_{k}(\omega )\leqslant
    \omega (\frac{\pi }{k}),\ \ \ k\in \mathbb{N},$$
formula (\ref{27.03.10-15:23:27}) is an asymptotic equality if and only if
\begin{equation}\label{23.12.10-00:16:44}
    \lim\limits_{t\to 0}{\frac{\omega (t)}{t}}=\infty.
\end{equation}
An example of moduli of continuity $\omega (t)$ which satisfy (\ref{23.12.10-00:16:44}), are
the functions:
    $$\omega (t)=t^\alpha ,\ \ \alpha \in(0,1),$$
    $$\omega (t)=\ln^\alpha (t+1),\ \ \alpha  \in(0,1),$$
    $$\omega (t)=
  \begin{cases}
    0, & t=0, \\
    t^\alpha \ln(\frac{1}{t}), & t\in(0,e^{-1/\alpha }],\ \ \ \alpha \in(0,1],\\
    \frac{1}{\alpha e}, & t\in[e^{-1/\alpha },\infty),
  \end{cases}
    $$
    $$\omega (t)=
  \begin{cases}
    0, & t=0, \\
    \ln^{-\alpha} (\frac{1}{t}), & t\in(0,e^{-(1+\alpha) }],\ \ \ \alpha \in(0,1],\\
    \frac{1}{(1+\alpha )^\alpha }, & t\in[e^{-(1+\alpha ) },\infty).
  \end{cases}
    $$

Putting $\omega (t)=t^\alpha ,$ $\alpha \in(0,1)$, in the hypothesis of Theorem
\hyperlink{15.04.11-16:12:48}{1} and taking into account that in this case the class
$H_\omega $ becomes the well-known H\"{o}lder class $H^\alpha ,$ we obtain the next
statement.

\textbf{Theorem 2.} \emph{Let $q\in(0,1),$ $\beta \in\mathbb{R}$,
 $n,p\in \mathbb{N},$
$p<n$ and $\alpha \in(0,1)$. Then the following asymptotic equality
\begin{equation}\label{27.04.10-10:27:17}
    \mathcal{E}(C^q_\beta H^\alpha  ;V_{n,p})=$$
    $$=
    \frac{q^{n-p+1}}{p(n-p+1)^\alpha }\Bigg(\frac{2^\alpha }{\pi ^2}
    {K}_{p,q}\int_{0}^{\pi /2}t^\alpha
    \sin t\,dt+\frac{O(1)}{(1-q)^{\delta (p)}(n-p+1)^{1-\alpha } }\Bigg),
    \end{equation}
is true as  $n-p\to\infty$, where ${K}_{p,q}$ and $\delta (p)$ are defined by
$(\ref{29.03.10-12:30:45})$ and $(\ref{29.03.10-12:31:28})$ respectively, and $O(1)$ is a
quantity uniformly bounded in $n,$ $p,$ $q$, $\alpha $ and $\beta .$}

We note that asymptotic equality (\ref{27.04.10-10:27:17}) for de la Vall\'{e}e Poussin sums
$V_{n,p}$ with bounded $p\in \mathbb{N}\setminus\{1\}$ is obtained for the first time.

Asymptotic behavior of the constant $K_{p,q}$ as $p\to\infty$ can be judged from the
estimate
    \begin{equation}\label{11.01.11-12:43:01}
    {K}_{p,q}=\frac{2\pi }{1-q^2}\Big(1+O(1)q^p\Big),
    \end{equation}
uniform in $p$ and $q$  (see \cite[p. 130]{Serdyuk V-P}). The substitution of
(\ref{11.01.11-12:43:01}) into equality (\ref{27.03.10-15:23:27}) enables us to obtain the
asymptotic estimate of the quantity $\mathcal{E}(C^q_\beta H_\omega ;V_{n,p})$ as
$n-p\to\infty$ in which the principal term coincides with the first summand in the
right-hand side of (\ref{7.05.10-12:08:22}).

In the general case $(p=1,2,\ldots,n)$, as follows from \cite[p. 215]{Savchuks and
Chaichenko}, the values of the constant $K_{p,q}$ can be expressed through the values of the
complete elliptic integral of the 1st kind $\textbf{K}(q^p)$ by means of the following
equation:
\begin{equation}\label{15.02.11-11:48:33}
    K_{p,q}=4\frac{1-q^{2p}}{1-q^2}\textbf{K}(q^p),\ \ p\in\mathbb{N},\ \ q\in(0,1).
\end{equation}
Using (\ref{15.02.11-11:48:33}), formulas (\ref{27.03.10-15:23:27}) and
(\ref{27.04.10-10:27:17}) have the next representations
    $$\mathcal{E}(C^q_\beta H_\omega ;V_{n,p})=
    \frac{q^{n-p+1}}{p}\Bigg(
    \frac{4}{\pi ^2}\frac{1-q^{2p}}{1-q^2}\textbf{K}(q^p)e_{n-p+1}(\omega )
    +\frac{O(1)\omega (\pi )}{(1-q)^{\delta (p)}(n-p+1)}\Bigg),\eqno{(9')}$$
    $$\mathcal{E}(C^q_\beta H^\alpha  ;V_{n,p})=$$
    $$=
    \frac{q^{n-p+1}}{p(n-p+1)^\alpha }\Bigg(\frac{2^{\alpha+2} }{\pi ^2}
    \frac{1-q^{2p}}{1-q^2}\textbf{K}(q^p)\int_{0}^{\pi /2}\!\!t^\alpha
    \sin t\,dt+\frac{O(1)}{(1-q)^{\delta (p)}(n-p+1)^{1-\alpha } }\Bigg).
    \eqno{(13')}$$

\hypertarget{15.04.11-16:24:36}{}

In addition to Theorem \hyperlink{15.04.11-16:12:48}{1} we present the sequent result.

\textbf{Theorem 3.} \emph{Let $q\in(0,1),$ $\beta \in\mathbb{R},$
 $n,p\in \mathbb{N},$
$p<n$ and let $\omega (t)$ be a convex upwards modulus of continuity. Then as
$n-p\to\infty$$:$
\begin{equation}\label{3.03.11-16:11:11}
    \mathcal{E}(C^q_\beta H_\omega ;V_{n,p}) =$$
    $$= \frac{q^{n-p+1}}{p}\Bigg(
    \frac{4J_{q,p}}{\pi ^2}\frac{1-q^p}{1-q}e_{n-p+1}(\omega )+\frac{O(1)}{(1-q)^{\delta (p)}
    (n-p+1)}\omega \bigg(\frac{1}{n-p+1}\bigg)\Bigg),
\end{equation}
where the two-sided estimate
\begin{equation}\label{3.03.11-16:17:10}
    \frac{1+q^{p}}{1+q}\mathbf{K}(q^p)\leqslant J_{q,p}\leqslant \mathbf{K}(q),
\end{equation}
holds for $J_{q,p}=J_{q,p}(n,\omega )$, and $e_{n-p+1}(\omega )$, $\delta (p)$ and $O(1)$
have the same meaning as in Theorem \hyperlink{15.04.11-16:12:48}{1}.}

For $p=1$, as appears from (\ref{3.03.11-16:17:10}), $J_{q,p}= J_{q,1}(n,\omega
)=\textbf{K}(q).$ In this case relation (\ref{3.03.11-16:11:11}) becomes the asymptotic
equality which is the special case of (\ref{11.05.10-12:13:35}).

\section{Proof of main result}

The proof of Theorem \hyperlink{15.04.11-16:12:48}{1} consists of three steps.

\emph{Step 1.} We single out a principal value of the quantity $\mathcal{E}( C^q_\beta
H_\omega;V_{n,p}).$

To this end we consider first the deviation
    \begin{equation}\label{22.03.10-12:33:40}
        \rho _{n,p}(f;x):=f(x)-V_{n,p}(f;x), \ \ f\in C^q_\beta H_\omega
    \end{equation}
and on the basis of equalities (\ref{6.05.10-14:46:39}) and (\ref{6.05.10-14:55:11}) we
write the representation
    \begin{equation}\label{29.03.10-13:16:27}
    \rho _{n,p}(f;x)=\frac{
    1}{\pi p}\int_{0}^{2\pi }\varphi (x+t)\sum\limits_{
    k=n-p}^{n-1}P_{q,\beta ,k+1}(t)\,dt,
    \end{equation}
in which $\varphi \in H_\omega ,$ and
$$
    P_{q,\beta ,m}(t)=\sum\limits_{j=m}^{\infty}q^j\cos\Big(jt+\frac{
    \beta \pi }{2}\Big),\ m\in \mathbb{N},\ q\in(0,1),\ \beta \in
    \mathbb{R}.
$$
Since
    $$\int_{0}^{2\pi }\sum\limits_{
    k=n-p}^{n-1}P_{q,\beta ,k+1}(t)\,dt=0$$
and, according to \cite[p. 118]{Stepanets_2001},
    $$P_{q,\beta ,m}(t)=q^mZ_q(t)\cos\Big(mt+\theta_q (t)+\frac{\beta \pi }{2}\Big),$$
where
    $$Z_{q}(t):=\frac{1}{\sqrt{1-2q\cos t+q^{2}}},$$
    $$\theta _{q}(t):=\text{arctg}\frac{q\sin t}{1-q\cos t},$$
it follows from (\ref{29.03.10-13:16:27}) that
    \begin{equation}\label{29.03.10-13:17:52}
    \rho _{n,p}(f;x)=\frac{
    1}{\pi p}\int_{0}^{2\pi }\big(\varphi (x+t)-\varphi (x)\big)Z_q(t)\sum\limits_{
    k=n-p+1}^{n}q^k\cos\Big(kt+\theta _q(t)+\frac{\beta \pi }{2}\Big)\,dt.
    \end{equation}
By virtue of formula (17) in \cite[p. 126]{Serdyuk V-P} the equality
    $$\sum\limits_{
    k=n-p+1}^{n}q^k\cos\Big(kt+\theta _q(t)+\frac{\beta \pi }{2}\Big)=$$
    $$=Z_q(t)q^{n-p+1}\bigg(\cos\Big((n-p+1)t+\frac{\beta \pi }{2}\Big)
    G_{p,q}(t)-$$
    \begin{equation}\label{29.03.10-13:39:58}
    -\sin\Big((n-p+1)t+\frac{\beta \pi }{2}\Big)
    H_{p,q}(t)\bigg)
    \end{equation}
holds, where
    $$G_{p,q}(t)=\cos 2\theta _q(t)-q^p\cos(pt+2\theta _q(t)),$$
    $$H_{p,q}(t)=\sin 2\theta _q(t)-q^p\sin(pt+2\theta _q(t)).$$
Representing the functions $G_{p,q}(t)$ and $H_{p,q}(t)$ in the form
    $$G_{p,q}(t)=\frac{\cos(2\theta _q(t)-\theta _{q^p}(pt))}{Z_{q^p}(pt)},\ \
    H_{p,q}(t)=\frac{\sin(2\theta _q(t)-\theta _{q^p}(pt))}{Z_{q^p}(pt)},$$
where
    $$Z_{q^p}(t):=\frac{1}{\sqrt{1-2q^p\cos t+q^{2p}}},\ \ \theta
    _{q^p}(t):=\text{arctg}\frac{q^p\sin t}{1-q^p\cos t},$$
from (\ref{29.03.10-13:39:58}) we get
\begin{equation}\label{29.03.10-13:55:16}
    \sum\limits_{
    k=n-p+1}^{n}q^k\cos\Big(kt+\theta _q(t)+\frac{\beta \pi }{2}\Big)=$$
    $$=q^{n-p+1}\frac{Z_q(t)}{Z_{q^p}(pt)}\cos\Big(
    (n-p+1)t+2\theta _q(t)-\theta _{q^p}(pt)+\frac{\beta \pi }{2}\Big).
\end{equation}
Relations (\ref{29.03.10-13:17:52}) and (\ref{29.03.10-13:55:16}) imply the following
equality
        \begin{equation}\label{22.03.10-12:48:22}
        \rho _{n,p}(f;x)=$$
        $$=\frac{q^{n-p+1}}{\pi p}\int_{0}^{
        2\pi }\big(\varphi (x+t)-\varphi (x)\big)\frac{Z^2_q(t)}{Z_{q^p}(pt)}\cos\Big(
        (n-p+1)t+2\theta _q(t)-\theta _{q^p}(pt)+\frac{\beta \pi
        }{2}\Big)\,dt.
        \end{equation}
The right-hand side of (\ref{22.03.10-12:48:22}) is 4-periodic in $\beta .$ Therefore, it
will be assumed below that $\beta \in[0,4).$

Since for any $\varphi \in H_\omega $ the function $\varphi _1(u)=\varphi (u+h),$ $h\in
\mathbb{R}$, also belongs to $H_\omega $,  then from (\ref{22.03.10-12:48:22}) we have
    \begin{equation}\label{22.03.10-12:49:30}
        \mathcal{E}(C^q_\beta H_\omega ;V_{n,p})=$$
        $$=\frac{q^{n-p+1}}{\pi p}
        \sup\limits_{\varphi \in H_\omega }{\bigg|\int_{0}^{
        2\pi }\!\!\Delta (\varphi ,t)\frac{Z^2_q(t)}{Z_{q^p}(pt)}\cos\Big(
        (n-p+1)t\!+2\theta _q(t)\!\!-\!\theta _{q^p}(pt)\!+\frac{\beta \pi
        }{2}\Big)\,dt\bigg|,}
    \end{equation}
where
    $$\Delta (\varphi ,t):=\varphi (t)-\varphi (0).$$
Denote by  $\mathcal{J}_{n,p,q,\beta }(\varphi )$ the integral in the right-hand side of
(\ref{22.03.10-12:49:30}), that is
\begin{equation}\label{22.03.10-12:56:41}
    \mathcal{J}_{n,p,q,\beta }(\varphi )=\int_{0}^{
        2\pi }\Delta (\varphi ,t)\frac{Z^2_q(t)}{Z_{q^p}(pt)}\cos\Big(
        (n-p+1)t\!+2\theta _q(t)\!\!-\!\theta _{q^p}(pt)\!+\frac{\beta \pi
        }{2}\Big)\,dt.
\end{equation}
Our further goal is to find an asymptotic estimation of the integral
$\mathcal{J}_{n,p,q,\beta }(\varphi )$ as \mbox{$n-p\to\infty.$} For this reason, without
loss of generality, we shall assume that the numbers $n$ and $p$ have been chosen such that
    \begin{equation}\label{22.03.10-13:20:36}
        n-p\geqslant \frac{6}{1-q}.
    \end{equation}
First we show that
    \begin{equation}\label{22.03.10-14:20:01}
    \mathcal{J}_{n,p,q,\beta }(\varphi )=$$
    $$=\int_{0}^{
        2\pi }\Delta (\varphi ,t)\frac{Z^4_q(t)}{Z_{q^p}(pt)Z_{q,n,p}^2(t)}\cos\Big(
        (n-p+1)t\!+2\theta _q(t)\!\!-\!\theta _{q^p}(pt)\!+\frac{\beta \pi
        }{2}\Big)\,dt+$$
        $$+O(1)\frac{\omega (\pi )}{(1-q)^{\delta (p)}(n-p+1)},
\end{equation}
where
\begin{equation}\label{30.11.10-12:35:06}
Z_{q,n,p}(t):=\frac{Z_q(t)}{
    \sqrt{\frac{n-p+1+2q(\cos t-q)Z_q^2(t)-pq^p(\cos pt-q^p)Z_{q^p}^2(pt)}{
    n-p+\alpha _q}}},
\end{equation}
    \begin{equation}\label{22.03.10-13:19:42}
    \alpha _q:=\Big[\frac{3q}{1-q}\Big]+2,
    \end{equation}
and $[ m]$ denotes the integral part of the number $m.$ By virtue of
(\ref{22.03.10-13:20:36}) and the obvious inequality
    \begin{equation}\label{7.12.10-15:06:27}
    pq^p|\cos pt-q^{p}|Z^2_{q^p}(pt)\leqslant \frac{pq^p}{1-q^p}
    \leqslant \frac{q}{1-q}, \ \ q\in(0,1),\ \ p\in\mathbb{N},\
    \ t\in \mathbb{R},
    \end{equation}
the quantity under the radical sign in (\ref{30.11.10-12:35:06}) is always positive. Set
    $$R_{n,p,q,\beta }(\varphi ):=\mathcal{J}_{n,p,q,\beta }(\varphi )-$$
    $$-\int_{0}^{
        2\pi }\Delta (\varphi ,t)\frac{Z^4_q(t)}{Z_{q^p}(pt)Z_{q,n,p}^2(t)}\cos\Big(
        (n-p+1)t\!+2\theta _q(t)\!\!-\!\theta _{q^p}(pt)\!+\frac{\beta \pi
        }{2}\Big)\,dt=$$
    \begin{equation}\label{11.01.11-12:50:03}
    =\int_{0}^{
        2\pi }\Delta (\varphi ,t)\frac{Z^2_q(t)}{Z_{q^p}(pt)}\bigg(1-\frac{
        Z_q^2(t)}{Z_{q,n,p}^2(t)}\bigg)\cos\Big(
        (n-p+1)t\!+2\theta _q(t)\!\!-\!\theta _{q^p}(pt)\!+\frac{\beta \pi
        }{2}\Big)\,dt.
    \end{equation}
To prove equality (\ref{22.03.10-14:20:01}) it suffices to show that the estimate
\begin{equation}\label{22.03.10-13:24:31}
    R_{n,p,q,\beta }(\varphi )=O(1)\frac{\omega (\pi )}{(1-q)^{\delta
    (p)}(n-p+1)}
\end{equation}
holds, where the quantity $\delta (p)$ is defined by (\ref{29.03.10-12:31:28}). Indeed,
considering the inequality $|\Delta (\varphi ,t)|\leqslant \omega (|t|),$ from
(\ref{11.01.11-12:50:03}) we get
\begin{equation}\label{22.03.10-13:28:52}
    \big|R_{n,p,q,\beta }(\varphi )\big|\leqslant \omega (2\pi )
    \int_{0}^{2\pi }\frac{Z^2_q(t)}{Z_{q^p}(pt)}\bigg|1-\frac{
        Z_q^2(t)}{Z_{q,n,p}^2(t)}\bigg|\,dt.
\end{equation}
After performing elementary transformations and taking into account
(\ref{7.12.10-15:06:27}), we find
    $$\bigg|1-\frac{
        Z_q^2(t)}{Z_{q,n,p}^2(t)}\bigg|=\bigg|\frac{2q(\cos t-q)Z_q^2(t)-
        pq^p(\cos pt-q^p)Z_{q^p}^2(pt)+1-\alpha _q}{
    n-p+\alpha _q}\bigg|< $$
    $$< \frac{1}{n-p+\alpha _q}\bigg(
    2q|\cos t-q|Z_{q}^2(t)+
    pq^p|\cos pt-q^p|Z_{q^p}^2(pt)+\frac{3q}{1-q}+1\bigg)\leqslant $$
    \begin{equation}\label{22.03.10-14:02:34}\leqslant \frac{1}{n-p+\alpha _q}\bigg(
    \frac{2q}{1-q}+\frac{
    q}{1-q}+\frac{3q}{1-q}+1\bigg)=
    \frac{O(1)}{(1-q)(n-p+1)}.
    \end{equation}
From the estimate
    $$\frac{Z^2_q(t)}{Z_{q^p}(pt)}=\frac{1}{\sqrt{
    {1-2q\cos t+q^2}}}\leqslant \frac{1}{1-q},\ \ t\in \mathbb{R},\ \ p=1,$$
and (\ref{11.01.11-12:43:01}) it follows that
\begin{equation}\label{22.03.10-14:08:31}
    \frac{Z^2_q(t)}{Z_{q^p}(pt)}=\frac{
    \sqrt{1-2q^p\cos pt+q^{2p}}}{1-2q\cos t+q^2}=\frac{O(1)}{(1-q
    )^{\delta (p)-1}},\ \ t\in \mathbb{R},\ \ p\in \mathbb{N}.
\end{equation}
Comparing (\ref{22.03.10-13:28:52})--(\ref{22.03.10-14:08:31}), we obtain
(\ref{22.03.10-13:24:31}), and with it estimate (\ref{22.03.10-14:20:01}).

Consider the function
    \begin{equation}\label{1.03.10-14:50:38}
        y_1(t):=t+\frac{1}{n-p+\alpha _q}\Big(2\theta _q(t)-\theta _{q^p}(pt)+(1-\alpha _q)t+\frac{
        \beta \pi }{2}\Big).
    \end{equation}
On the strength of the fact that
    $(\theta _{q^p}(pt))'=pq^p(\cos pt-q^p)Z_{q^p}^2(pt),$ the equality
\begin{equation}\label{1.03.10-14:54:19}
    y_1'(t)=1+\frac{1}{n-p+\alpha _q}\Big(2q(\cos t-q)Z_q^2(t)-pq^p(\cos pt-q^p)Z_{q^p}^2
    (pt)+1-\alpha _q\Big)=$$
    $$=\frac{Z_q^2(t)}{Z_{q,n,p}^2(t)},
\end{equation}
holds, where $Z_{q,n,p}(t)$ is defined by (\ref{30.11.10-12:35:06}). In view of
(\ref{7.12.10-15:06:27})
\begin{equation}\label{11.01.11-15:00:41}
    -\frac{5q+1}{1-q}\leqslant 2q(\cos t-q)Z_q^2(t)-pq^p(\cos pt-q^p)Z_{q^p}^2
    (pt)+1-\alpha _q<0.
\end{equation}
Thus, we see from (\ref{22.03.10-13:20:36}) and (\ref{1.03.10-14:54:19}) that the next
two-sided estimate
\begin{equation}\label{22.03.10-14:54:30}
    \frac{1}{3}<y_1'(t)<1
\end{equation}
holds. So $y_1$ has the inverse function $y(t)=y_1^{-1}(t),$ whose derivative $y'$ by
(\ref{1.03.10-14:54:19}) satisfies the relation
    \begin{equation}\label{1.04.10-14:54:26}
    y'(t)=\frac{1}{y_1'(y(t))}=\frac{Z_{q,n,p}^2(y(t))}{Z_q^2(y(t))}.
    \end{equation}
Making the change of variables $t=y(\tau )$ in (\ref{22.03.10-14:20:01}), we obtain by
(\ref{1.03.10-14:54:19}) the relation
\begin{equation}\label{1.03.10-15:53:26}
    \mathcal{J}_{n,p,q,\beta }(\varphi )=$$
    $$=\int_{y_1(0)}^{y_1(2\pi )}\Delta(\varphi ,y(\tau ))
    \frac{Z_q^4(y(\tau ))}{Z_{q^p}(py(\tau ))Z_{q,n,p}^2(y(\tau ))}
    \frac{Z_{q,n,p}^2(y(\tau ))}{Z_q^2(y(\tau ))}\cos\big((n-p+\alpha _q)\tau\big) \,d\tau+$$
        $$+O(1)\frac{\omega (\pi )}{(1-q)^{\delta (p)}(n-p+1)}=$$
    $$=\int_{y_1(0)}^{y_1(2\pi )}\Delta
    (\varphi ,y(\tau ))
    \frac{Z_q^2(y(\tau ))}{Z_{q^p}(py(\tau ))}
    \cos\big((n-p+\alpha _q)\tau\big) d\tau+$$
    $$+O(1)\frac{\omega (\pi )}{(1-q)^{\delta (p)}(n-p+1)}.
\end{equation}
We set
    \begin{equation}\label{27.01.11-10:33:05}
    x_k:=\frac{k\pi }{n-p+\alpha _q},\ \ \tau _k:=x_k+\frac{\pi }{2(n-p+\alpha _q)},\ \ k\in\mathbb{N}
    \end{equation}
and
    $$l_n(\tau )=
  \begin{cases}
    \frac{Z_{q}^2(y(\tau _k))}{Z_{q^p}(py(\tau _k))}, & \tau \in[x_k,x_{k+1}],\
    k=2,3,\ldots,k_0-2,k_0-1, \\
    0, & \tau \in[y_1(0),x_2)\cup(x_{k_0},y_1(2\pi )],
  \end{cases}
    $$
where $k_0$ is an index such that $\tau _{k_0}$ is the nearest to the left of $y_1(2\pi )$
root of the function $\cos\big((n-p+\alpha _q)\tau\big).$ Using this notations, from
(\ref{1.03.10-15:53:26}) we get
\begin{equation}\label{22.03.10-16:00:30}
    \mathcal{J}_{n,p,q,\beta }(\varphi )=\int_{
    x_2}^{x_{k_0}}\Delta (\varphi ,y(\tau ))l_n(\tau )\cos\big((n-p+\alpha _q)\tau \big)\,d\tau
    +R_{n,p,q }^{(1)}(\varphi )+R_{n,p,q }^{(2)}(\varphi )+$$
    $$+O(1)\frac{\omega (\pi )}{(1-q)^{\delta (p)}(n-p+1)},
\end{equation}
where
    $$R_{n,p,q }^{(1)}(\varphi ):=$$
    $$=\sum\limits_{k=2}^{k_0-1}\int_{
    x_k}^{x_{k+1}}\Delta (\varphi ,y(\tau ))\bigg(
    \frac{Z_{q}^2(y(\tau ))}{Z_{q^p}(py(\tau ))}-
    \frac{Z_{q}^2(y(\tau _k))}{Z_{q^p}(py(\tau _k))}\bigg)\!\!\cos\big((n-p+\alpha _q)\tau\big)
    \,d\tau,$$
    $$R_{n,p,q }^{(2)}(\varphi ):=\bigg(\int_{y_1(0)}^{x_2}+\int_{x_{k_0}}^{y_1(2\pi )}\bigg)
    \Delta (\varphi ,y(\tau ))
    \frac{Z_{q}^2(y(\tau ))}{Z_{q^p}(py(\tau ))}\cos\big((n-p+\alpha _q)\tau\big)
    \,d\tau.$$
Since $y(t)$ is an increasing function (see (\ref{22.03.10-14:54:30}),
(\ref{1.04.10-14:54:26})), then
    $$y(x_{k+1})\leqslant y(x_{k_0})<y(y_1(2\pi ))=2\pi ,\ \ \ k=\overline{2,k_{0}-1}$$
and so for $R_{n,p,q }^{(1)}(\varphi )$ the following trivial estimate
\begin{equation}\label{1.04.10-14:50:54}
    \big|R_{n,p,q }^{(1)}(\varphi )\big|\leqslant \sum\limits_{k=2}^{k_0-1}
    \omega (y(x_{k+1}))\int_{x_k}^{x_{k+1}}\bigg|
    \frac{Z_{q}^2(y(\tau ))}{Z_{q^p}(py(\tau ))}-
    \frac{Z_{q}^2(y(\tau _k))}{Z_{q^p}(py(\tau _k))}\bigg|\,d\tau<$$
    $$<\omega (2\pi )\sum\limits_{k=2}^{k_0-1}
    \int_{x_k}^{x_{k+1}}\bigg|
    \frac{Z_{q}^2(y(\tau ))}{Z_{q^p}(py(\tau ))}-
    \frac{Z_{q}^2(y(\tau _k))}{Z_{q^p}(py(\tau _k))}\bigg|\,d\tau
\end{equation}
holds.

We will show that
\begin{equation}\label{6.04.10-11:08:09}
    \sum\limits_{k=2}^{k_0-1}
    \int_{x_k}^{x_{k+1}}\bigg|
    \frac{Z_{q}^2(y(\tau ))}{Z_{q^p}(py(\tau ))}-
    \frac{Z_{q}^2(y(\tau _k))}{Z_{q^p}(py(\tau _k))}\bigg|\,d\tau
    =\frac{O(1)q}{(1-q)^{\delta (p)}(n-p+1)}.
\end{equation}
For this purpose, we consider the derivative
\begin{equation}\label{22.03.10-16:21:35}
    \bigg(\frac{Z_q^2(t)}{Z_{q^p}(pt)}\bigg)'=
    -2q\sin t \frac{Z_q^4(t)}{Z_{q^p}(pt)}+q^pp\sin pt\, Z_{q^p}(pt)Z_q^2(t)=:J_{q,p}^{(1)}(t)+J_{q,p}^{(2)}(t)
\end{equation}
and estimate separately the summands of the right-hand side of (\ref{22.03.10-16:21:35}).
Taking into account (\ref{22.03.10-14:08:31}) and the estimate
\begin{equation}\label{22.03.10-17:47:21}
    q|\sin t|Z_q^2(t)=\frac{q|\sin t|}{1-2q\cos t+q^2}<\sum\limits_{k=1}^{\infty}q^k=\frac{q}{1-q},
    \end{equation}
 we obtain for $J_{q,p}^{(1)}(t)$
    \begin{equation}\label{22.03.10-16:48:03}
    \big|J_{q,p}^{(1)}(t)\big|\leqslant 2q|\sin t|Z_q^2(t)
    \frac{Z_q^2(t)}{Z_{q^p}(pt)}=O(1)\frac{q}{(1-q)^{\delta (p)}}.
    \end{equation}
If $p=1$, then in view of (\ref{22.03.10-17:47:21}) we have for $J_{q,p}^{(2)}(t)$
    $$
    \big|J_{q,1}^{(2)}(t)\big|<\frac{q}{(1-q)^2}.
    $$
For any $p=2,3,\ldots,$ one easily shows that
    $$\big|J_{q,p}^{(2)}(t)\big|<\frac{q^pp}{(1-q)^2(1-q^p)}<\frac{q}{(1-q)^3}.$$
Therefore we finally obtain
\begin{equation}\label{22.03.10-17:52:11}
    \big|J_{q,p}^{(2)}(t)\big|<\frac{q}{(1-q)^{\delta (p)}},\ \ p\in \mathbb{N}.
\end{equation}
Combining (\ref{22.03.10-16:21:35}), (\ref{22.03.10-16:48:03}) and
(\ref{22.03.10-17:52:11}), we arrive at the estimate
    $$\bigg|\bigg(\frac{Z_q^2(t)}{Z_{q^p}(pt)}\bigg)'\bigg|=O(1)\frac{
    q}{(1-q)^{\delta (p)}},\ \ t\in [0,2\pi ].$$
Since by (\ref{22.03.10-14:54:30}) and (\ref{1.04.10-14:54:26}) $|y'(t)|<3$, applying
Lagrange's theorem on finite increments, we find
\begin{equation}\label{22.03.10-18:03:54}
    \bigg|
    \frac{Z_{q}^2(y(\tau ))}{Z_{q^p}(py(\tau ))}-
    \frac{Z_{q}^2(y(\tau _k))}{Z_{q^p}(py(\tau _k))}\bigg|=$$
    $$=
    O(1)\frac{
    q}{(1-q)^{\delta (p)}}|y(\tau )-y(\tau _k)|=O(1)\frac{
    q}{(1-q)^{\delta (p)}}|\tau -\tau _k|=$$
    $$=O(1)\frac{
    q}{(1-q)^{\delta (p)}(n-p+\alpha _q)}=$$
    $$=O(1)\frac{
    q}{(1-q)^{\delta (p)}(n-p+1)},\ \ \tau \in[x_k,x_{k+1}],\ \
    k=\overline{2,k_0-1}.
\end{equation}
It results from (\ref{22.03.10-18:03:54}) that
\begin{equation}\label{6.04.10-12:34:50}
    \sum\limits_{k=2}^{k_0-1}
    \int_{x_k}^{x_{k+1}}\bigg|
    \frac{Z_{q}^2(y(\tau ))}{Z_{q^p}(py(\tau ))}-
    \frac{Z_{q}^2(y(\tau _k))}{Z_{q^p}(py(\tau _k))}\bigg|\,d\tau
    =O(1)\frac{qx_{k_0}}{(1-q)^{\delta (p)}(n-p+1)}.
\end{equation}
But since
$$
    x_{k_0}<y_1(2\pi )=\int_{0}^{2\pi }y_1'(t)\,dt+y_1(0)=$$
    $$=\int_{0}^{2\pi }y_1'(t)\,dt+\frac{\beta \pi }{2(n-p+\alpha _q)}<2\pi +\frac{
    \beta \pi }{2}< 4\pi ,
$$
(\ref{6.04.10-11:08:09}) follows from (\ref{6.04.10-12:34:50}). Estimates
(\ref{1.04.10-14:50:54}) and (\ref{6.04.10-11:08:09}) imply
\begin{equation}\label{22.03.10-18:16:13}
    \big|R_{n,p,q }^{(1)}(\varphi )\big|=O(1)\frac{\omega (\pi )q}{(1-q)^{\delta (p)}(n-p+1)}.
\end{equation}
Further, considering that
    \begin{equation}\label{24.03.10-11:55:27}
    x_2-y_1(0)\leqslant \frac{2\pi }{n-p+\alpha _q}<\frac{2\pi }{n-p+1},
    \end{equation}
    \begin{equation}\label{24.03.10-11:55:44}
    y_1(2\pi )-x_{k_0}\leqslant \tau _{k_0+1}-x_{k_0}=\frac{3\pi }{2(n-p+\alpha _q)}<
    \frac{3\pi }{2(n-p+1)}
    \end{equation}
and using (\ref{22.03.10-14:08:31}), we find
\begin{equation}\label{22.03.10-18:53:03}
    \big|R_{n,p,q }^{(2)}(\varphi )\big|
    =O(1)\frac{\omega (\pi )}{(1-q)^{\delta (p)-1}(n-p+1)}.
\end{equation}
From (\ref{22.03.10-16:00:30}), (\ref{22.03.10-18:16:13}) and (\ref{22.03.10-18:53:03}) we
obtain
\begin{equation}\label{22.03.10-19:05:11}
    \mathcal{J}_{n,p,q,\beta }(\varphi )=\int_{
    x_2}^{x_{k_0}}\Delta (\varphi ,y(\tau ))l_n(\tau )
   \cos \big((n-p+\alpha _q)\tau\big) \,d\tau+$$
    $$+O(1)\frac{\omega (\pi )}{(1-q)^{\delta (p)}(n-p+1)},\ \ \varphi \in H_\omega ,
    \ \ n-p\to\infty,
\end{equation}
where $O(1)$ is quantity uniformly bounded relative to all parameters under consideration.

Comparing (\ref{22.03.10-12:49:30}), (\ref{22.03.10-12:56:41}) and (\ref{22.03.10-19:05:11})
we conclude that
\begin{equation}\label{23.03.10-19:08:49}
        \mathcal{E}(C^q_\beta H_\omega ;V_{n,p})=\frac{q^{n-p+1}}{\pi p}
        \bigg(\sup\limits_{\varphi \in H_\omega }{\big|I_{n,p,q,\beta }(\varphi )
        \big|+\frac{O(1)\omega (\pi )}{(1-q)^{\delta (p)}(n-p+1)}\bigg),}
    \end{equation}
in which
    \begin{equation}\label{27.03.10-13:59:05}
    I_{n,p,q,\beta }(\varphi ):=\int_{
    x_2}^{x_{k_0}}\Delta (\varphi ,y(\tau ))l_n(\tau )
    \cos\big((n-p+\alpha _q)\tau\big) \,d\tau=$$
    $$=\sum\limits_{k=2}^{k_0-1}
    \frac{Z_{q}^2(y(\tau _k))}{Z_{q^p}(py(\tau _k))}\int_{x_k}^{x_{k+1}}
    \Delta (\varphi ,y(\tau ))\cos\big((n-p+\alpha _q)\tau\big)\,d\tau.
    \end{equation}

\emph{Step 2.} Using formula (\ref{23.03.10-19:08:49}) we find an upper bound for
$\mathcal{E}(C^q_\beta H_\omega ;V_{n,p}).$

With this goal, dividing each integral
    $$\int_{x_k}^{x_{k+1}}\Delta (\varphi ,y(\tau ))\cos\big((n-p+\alpha _q)\tau\big) \,d\tau,\ \ k=\overline{2,k_0-1} $$
into two integrals over $(x_k, \tau _k)$ and $(\tau _k, x_{k+1})$, and setting $z=2\tau
_k-\tau $ in the last integral, we obtain
    $$\bigg|\int_{x_k}^{x_{k+1}}
    \Delta (\varphi ,y(\tau ))\cos\big((n-p+\alpha _q)\tau\big) \,d\tau \bigg|=$$
    \begin{equation}\label{20.07.10-02:43:22}
    =\bigg|\int_{x_k}^{\tau _{k}}\big(\varphi (y(\tau ))-
    \varphi (y(2\tau _k-\tau ))\big)
    \cos\big((n-p+\alpha _q)\tau\big) \,d\tau \bigg|,\ \ k=\overline{2,k_0-1}.
    \end{equation}
We choose $c_k\in[x_k,x_{k+1}]$ such that
    $$y'(c_k)=\max\limits_{\tau \in[x_k,x_{k+1}]}y'(\tau ).$$
Then for any $\tau \in[x_k,\tau _{k}]$
    $$y(2\tau _k-\tau )-y(\tau )=\int_{
    \tau }^{2\tau_k -\tau }y'(x)\,dx\leqslant y'(c_k)(2\tau _k-2\tau )$$
and consequently
\begin{equation}\label{20.07.10-02:38:17}
    \big|\varphi (y(\tau ))-
    \varphi (y(2\tau _k-\tau ))\big|\leqslant $$
    $$\leqslant \omega \big(y(2\tau _k-\tau )-y(\tau )\big)\leqslant
    \omega (2y'(c_k)(\tau _k-\tau )),\ \ k=\overline{2,k_0-1}.
\end{equation}
From (\ref{20.07.10-02:43:22}), in view of (\ref{20.07.10-02:38:17}),  we find that
    $$\bigg|\int_{x_k}^{x_{k+1}}
    \Delta (\varphi ,y(\tau ))\cos\big((n-p+\alpha _q)\tau\big) \,d\tau \bigg|\leqslant $$
    $$\leqslant
    \int_{x_k}^{\tau _{k}}\omega \big(2y'(c_k )(\tau _k-\tau )\big)
    \big|\cos\big((n-p+\alpha _q)\tau\big) \big|\,d\tau\leqslant $$
    \begin{equation}\label{17.05.10-14:42:37}
    \leqslant \frac{1}{n-p+\alpha _q}\int_{0}^{\pi /2}\omega \bigg(\frac{2
    y'(c_k)t}{n-p+\alpha _q}\bigg)\sin t\,dt<$$
    $$<\frac{1}{n-p+\alpha _q}\int_{0}^{\pi /2}\omega \bigg(\frac{2
    y'(c_k)t}{n-p+1}\bigg)
    \sin t\,dt=$$
    $$=\frac{1}{n-p+\alpha _q}\int_{0}^{\pi /2}\omega \bigg(\frac{2
    t}{n-p+1}\bigg)
    \sin t\,dt+$$
    $$+\frac{O(1)}{n-p+\alpha _q}\max\limits_{t\in(0,\pi /2]}{\Big|\omega \bigg(\frac{
    2y'(c_k)t}{n-p+1}\bigg)-\omega \bigg(\frac{
    2t}{n-p+1}\bigg)\Big|},\ \ k=\overline{2,k_0-1}.
    \end{equation}
Because for any convex upwards modulus of continuity $\omega $
    \begin{equation}\label{10.12.10-14:56:23}
    \omega (b)-\omega (a)\leqslant \omega (a)\frac{b-a}{a},\ \ 0<a<b,
    \end{equation}
then taking into consideration that by (\ref{22.03.10-14:54:30}) and
(\ref{1.04.10-14:54:26}) $y'(c_k)>1$, we have
    $$\max\limits_{t\in(0,\pi /2]}{\Big|\omega \bigg(\frac{
    2y'(c_k)t}{n-p+1}\bigg)-\omega \bigg(\frac{
    2t}{n-p+1}\bigg)\Big|}\leqslant$$
    \begin{equation}\label{8.12.10-12:48:53}
    \leqslant  \omega \bigg(\frac{\pi }{n-p+1}\bigg)(y'(c_k)-1),\ \ k=\overline{2,k_0-1}.
    \end{equation}
By virtue of (\ref{1.04.10-14:54:26}),
    \begin{equation}\label{8.12.10-11:53:26}y'(c_k)-1=\frac{Z^2_{q,n,p}(y(c_k))}{Z^2_q(y(c_k))}-1=\frac{
    \alpha _q-1-\lambda _{p,q}(c_k)}{
    {n-p+1+\lambda _{p,q}(c_k)}},
    \end{equation}
where
    $$\lambda _{p,q}(c_k)=2q(\cos y(c_k)-q)Z_q^2(y(c_k))-pq^p(\cos py(c_k)-q^p)Z_{q^p}^2(py(c_k)).$$
In view of (\ref{7.12.10-15:06:27}),
    $$|\lambda _{p,q}(c_k)|\leqslant \frac{
    3q}{1-q}.$$
Hence it follows from (\ref{8.12.10-11:53:26}), taking into account
(\ref{22.03.10-13:20:36}), that
   \begin{equation}\label{8.12.10-12:49:29}
   y'(c_k)-1\leqslant \frac{\alpha _q-1+\frac{3q}{1-q}}{n-p+1-\frac{3q}{1-q}}
    <\frac{6}{(1-q)(n-p+1-\frac{3q}{1-q})}<\frac{12}{(1-q)(n-p+1)}.
   \end{equation}
Comparing (\ref{17.05.10-14:42:37}), (\ref{8.12.10-12:48:53}) and (\ref{8.12.10-12:49:29}),
we obtain
\begin{equation}\label{23.03.10-22:49:41}
    \bigg|\int_{x_k}^{x_{k+1}}
    \Delta (\varphi ,y(\tau ))\cos(n-p+\alpha _q)\tau \,d\tau \bigg|\leqslant $$
    $$\leqslant
    \frac{1}{n-p+\alpha _q}\Bigg(\int_{0}^{\pi /2}\omega \bigg(\frac{2t}{n-p+1}\bigg)
    \sin t\,dt+$$
    $$+\frac{O(1)}{(1-q)(n-p+1)}
    \omega\bigg(\frac{1}{n-p+1}\bigg)\Bigg),\ \ k=\overline{2,k_0-1}.
\end{equation}
Applying to each integral in (\ref{27.03.10-13:59:05}) estimate (\ref{23.03.10-22:49:41}),
we have
\begin{equation}\label{24.03.10-15:14:51}
    |I_{n,p,q,\beta }(\varphi )|\leqslant
    \frac{1}{n-p+\alpha _q}\sum\limits_{k=2}^{k_0-1}\frac{Z_q^2(y(\tau _k))}{Z_{q^p}(py(\tau _k))}\int_{0}^{\pi /2}\omega \bigg(\frac{2t}{n-p+1}\bigg)
    \sin t\,dt+$$
    $$+\frac{O(1)}{(1-q)^{\delta (p)}(n-p+1)}
    \omega\bigg(\frac{1}{n-p+1}\bigg).
\end{equation}
Let us show that
\begin{equation}\label{24.03.10-15:26:28}
    \frac{\pi }{n-p+\alpha _q}\sum\limits_{k=2}^{k_0-1}\frac{Z_q^2(y(\tau _k))}{Z_{q^p}(py(\tau
    _k))}=\int_{0}^{2\pi }\frac{Z_q^2(t)}{Z_{q^p}(pt)}\,dt+\frac{
    O(1)}{(1-q)^{\delta (p)}(n-p+1)}.
\end{equation}
For this, we represent the left-hand side of (\ref{24.03.10-15:26:28}) as
    \begin{equation}\label{12.01.11-14:37:01}
    \frac{\pi }{n-p+\alpha _q}\sum\limits_{k=2}^{k_0-1}\frac{Z_q^2(y(\tau _k))}{Z_{q^p}(py(\tau
    _k))}=\int_{
    y_1(0)}^{y_1(2\pi )}\frac{Z_q^2(y(\tau ))}{Z_{q^p}(py(\tau
    ))}\,d\tau +R_{n,p,q}^{(1)}+R_{n,p,q}^{(2)},
    \end{equation}
where
    $$R_{n,p,q}^{(1)}:=-\bigg(\int_{y_1(0)}^{x_2}+
    \int_{x_{k_0}}^{y_1(2\pi )}\bigg)\frac{Z_q^2(y(\tau ))}{Z_{q^p}(py(\tau
    ))}\,d\tau,$$
    $$R_{n,p,q}^{(2)}:=\sum\limits_{k=2}^{k_0-1}\int_{
    x_k}^{x_{k+1}}\bigg(
    \frac{Z_{q}^2(y(\tau _k))}{Z_{q^p}(py(\tau _k))}-\frac{Z_{q}^2(
    y(\tau ))}{Z_{q^p}(py(\tau ))}\bigg)\,d\tau.$$
By virtue of (\ref{22.03.10-14:08:31}), (\ref{24.03.10-11:55:27}) and
(\ref{24.03.10-11:55:44})
    \begin{equation}\label{12.01.11-14:35:16}
    R_{n,p,q}^{(1)}=\frac{O(1)}{(1-q)^{\delta (p)-1}(n-p+1)},
    \end{equation}
and by (\ref{6.04.10-11:08:09})
    \begin{equation}\label{12.01.11-14:35:29}
    R_{n,p,q}^{(2)}=\frac{O(1)q}{(1-q)^{\delta (p)}(n-p+1)}.
    \end{equation}
Combining (\ref{12.01.11-14:37:01})--(\ref{12.01.11-14:35:29}), we can write
    $$\frac{\pi }{n-p+\alpha _q}\sum\limits_{k=2}^{k_0-1}\frac{Z_q^2(y(\tau _k))}{Z_{q^p}(py(\tau
    _k))}=\int_{
    y_1(0)}^{y_1(2\pi )}\frac{Z_q^2(y(\tau ))}{Z_{q^p}(py(\tau
    ))}\,d\tau +\frac{O(1)}{(1-q)^{\delta (p)}(n-p+1)}=$$
    \begin{equation}\label{8.12.10-15:12:30}
    =
    \int_{
    0}^{2\pi }\frac{Z_q^2(t)}{Z_{q^p}(pt)}\,dt+
    \int_{
    0}^{2\pi }\frac{Z_q^2(t)}{Z_{q^p}(pt)}\big(y'_1(t)-1\big)\,dt
    +\frac{O(1)}{(1-q)^{\delta (p)}(n-p+1)}.
    \end{equation}
But in view of (\ref{1.03.10-14:54:19}) and (\ref{11.01.11-15:00:41})
    $$|y'_1(t)-1|<\frac{6}{(1-q)(n-p+1)}.$$
Thus, in consideration of (\ref{22.03.10-14:08:31}), we obtain from (\ref{8.12.10-15:12:30})
equality (\ref{24.03.10-15:26:28}).

Estimates (\ref{23.03.10-19:08:49}), (\ref{24.03.10-15:14:51}) and (\ref{24.03.10-15:26:28})
imply that as $n-p\to\infty$
\begin{equation}\label{24.03.10-16:07:09}
    \mathcal{E}(C^q_\beta H_\omega ;V_{n,p})\leqslant
    \frac{q^{n-p+1}}{\pi ^2p}
    {K}_{p,q}\int_{0}^{\pi /2}\omega \bigg(\frac{2t}{n-p+1}\bigg)
    \sin t\,dt+$$
    $$+O(1)\frac{q^{n-p+1}\omega (\pi )}{(1-q)^{\delta (p)}p(n-p+1)},
\end{equation}
where
\begin{equation}\label{27.01.11-12:02:06}
{K}_{p,q}=\int_{
    0}^{2\pi }\frac{Z_q^2(t)}{Z_{q^p}(pt)}\,dt,
\end{equation}
and $O(1)$ is a quantity uniformly bounded in $n,$ $p,$ $q$, $\omega $ and $\beta .$

\emph{Step 3.} We now show that (\ref{24.03.10-16:07:09}) is an equality. For this, in view
of (\ref{23.03.10-19:08:49}), it is sufficient to show that there exists a function $\varphi
^*\in H_\omega $ such that
    \begin{equation}\label{25.01.11-15:00:49}
       I_{n,p,q,\beta }(\varphi^* )=\frac{{K}_{p,q}}{\pi}
    \int_{0}^{\pi /2}\omega \bigg(\frac{2t}{n-p+1}\bigg)
    \sin t\,dt+\frac{O(1)\omega (\pi )}{(1-q)^{\delta (p)}(n-p+1)},
    \end{equation}
where $I_{n,p,q,\beta }(\varphi^* )$ is defined by (\ref{27.03.10-13:59:05}). To this end,
we set
    $$\varphi _i(t):=
  \begin{cases}
    \frac{1}{2}\omega \big(2y_1(t)-2\tau _i\big),  & t\in[y(\tau _i),y(x_{i+1})], \\
    \frac{1}{2}\omega \big(2\tau _{i+1}-2y_1(t
    )\big), & t\in[y(x_{i+1}),y(\tau _{i+1})],\ \ i=\overline{s,k_0-1},
  \end{cases} $$
  $$s=
  \begin{cases}
    2, &  \text{if } k_0 \text{ is odd}, \\
    3, & \text{if } k_0 \text{ is even},
  \end{cases}
  $$
    $$\widetilde{\varphi }(t):=(-1)^{i+1}\varphi _i(t),\ \ t\in
    [y(\tau _i),y(\tau _{i+1})],\ \ i=\overline{s,k_0-1}.$$
Since $\tau _{k_0}\leqslant y_1(2\pi )$ and by (\ref{27.01.11-10:33:05}),
(\ref{1.03.10-14:50:38}) and $\beta \in[0,4),$ the inequality $\tau_s>y_1(0)$ holds, it
follows that
  $y(\tau _{k_0})\leqslant 2\pi $ and $y(\tau _s)>0$.
Consider the function
    \begin{equation}\label{27.01.11-11:33:54}
    \varphi ^*(t):=
  \begin{cases}
    \widetilde{\varphi }(t), & t\in[y(\tau _s),y(\tau _{k_0})], \\
    0, & t\in[0,2\pi ]\setminus[y(\tau _s),y(\tau _{k_0})],
  \end{cases}\ \ \varphi ^*(t)=\varphi ^*(t+2\pi ).
    \end{equation}
We show that, if (\ref{22.03.10-13:20:36}) holds, then $\varphi ^*\in H_\omega .$ The
construction of $\varphi ^*$ shows that it suffices to establish the inequality
    $$|\varphi ^*(t')-\varphi ^*(t'')|\leqslant \omega (t''-t'),$$
where $t'\in[y(x_i),y(\tau _i)]$ and $t''\in[y(\tau _i),y(x_{i+1})],$
$i=\overline{s+1,k_{0}-1}.$

In view of the convexity (upwards) of the modulus of continuity
    $$\frac{1}{2}\big(
    \omega (t_1)+\omega (t_2)\big)\leqslant \omega
    \bigg(\frac{t_1+t_2}{2}\bigg).$$
Then, considering that by (\ref{22.03.10-14:54:30}) $y_1'\in(\frac{1}{3},1)$, we get
    $$|\varphi ^*(t')-\varphi ^*(t'')|=|\widetilde{\varphi }(t')-
    \widetilde{\varphi }(t'')|=\varphi _{i-1}(t')+\varphi _i(t'')=$$
    $$=\frac{1}{2}\big(\omega (
    2\tau _i-2y_1(t'))+\omega (2y_1(t'')-2\tau _i)\big)\leqslant $$
    $$\leqslant \omega (y_1(t'')-y_1(t'))\leqslant \omega (t''-t'),\ \ i=\overline{s+1,k_0-1}.$$
Let us now verify that $\varphi ^*(t)$ is the desired function, which means that $\varphi
^*(t)$ satisfies (\ref{25.01.11-15:00:49}). Since by (\ref{27.01.11-11:33:54})
    $$\varphi^*(y(\tau ))=
  \begin{cases}
    \frac{(-1)^{i+1}}{2}\omega (2\tau -2\tau _i), & \tau \in[\tau _i,x_{i+1}], \\
    \frac{(-1)^{i+1}}{2}\omega (2\tau _{i+1}-2\tau ), &
    \tau \in[x_{i+1},\tau _{i+1}],
  \end{cases}\ \ \ i=\overline{s,k_0-1},
    $$
it follows that
    $$\int_{x_k}^{x_{k+1}}\Delta (\varphi ^*,y(\tau ))\cos\big((n-p+\alpha _q)\tau\big)
    \,d\tau=$$
    $$=\frac{(-1)^k}{2}\bigg(\int_{
    x_k}^{\tau _k}\omega (2\tau _k-2\tau )
    \cos\big((n-p+\alpha _q)\tau\big) \,d\tau -$$
    $$-\int_{\tau _k}^{x_{k+1}}\omega (2\tau -2\tau _k)
    \cos\big((n-p+\alpha _q)\tau\big) \,d\tau \bigg)=$$
    $$=\int_{0}^{\pi /2(n-p+\alpha _q)}\omega (2t)\sin\big((n-p+\alpha _q
    )\tau\big) \,d\tau =$$
    \begin{equation}\label{23.04.10-14:38:52}
    =\frac{1}{n-p+\alpha _q}\int_{0}^{\pi /2}\omega \bigg(
    \frac{2t}{n-p+\alpha _q}\bigg)\sin t\,dt,\ \ k=\overline{s+1,k_0-1}.
    \end{equation}
By (\ref{23.04.10-14:38:52}) and (\ref{22.03.10-14:08:31}), we obtain
    $$I_{n,p,q,\beta }(\varphi^* )=\sum\limits_{k=2}^{k_0-1}
    \frac{Z_{q}^2(y(\tau _k))}{Z_{q^p}(py(\tau _k))}\int_{x_k}^{x_{k+1}}
    \Delta (\varphi^* ,y(\tau ))\cos\big((n-p+\alpha
    _q)\tau\big)\,d\tau=$$
    $$
    =
    \frac{
    1 }{n-p+\alpha _q}\sum\limits_{k=2}^{k_0-1}\frac{Z_q^2(y(\tau _k))}{Z_{q^p}(py(\tau
    _k))}\int_{0}^{\pi /2}\omega \bigg(
    \frac{2t}{n-p+\alpha _q}\bigg)\sin t\,dt+$$
    \begin{equation}\label{27.01.11-12:38:42}
    +\frac{O(1)}{(n-p+1)(1-q)^{\delta (p)-1}}\omega \bigg(
    \frac{1}{n-p+1}\bigg).
    \end{equation}
From (\ref{27.01.11-12:38:42}), in view of  (\ref{24.03.10-15:26:28}) and
(\ref{27.01.11-12:02:06}), we find
    $$I_{n,p,q,\beta }(\varphi^* )=\frac{{K}_{p,q}}{\pi }\int_{0}^{\pi /2}\omega \bigg(
    \frac{2t}{n-p+\alpha _q}\bigg)\sin t\,dt+$$
    \begin{equation}\label{27.01.11-12:38:05}
    +\frac{O(1)}{(n-p+1)(1-q)^{\delta (p)}}\omega \bigg(
    \frac{1}{n-p+1}\bigg).
    \end{equation}
Based on (\ref{10.12.10-14:56:23}) and (\ref{22.03.10-13:19:42}), we can write
    $$\max\limits_{t\in(0,\pi /2]}{\bigg|\omega \bigg(\frac{2t}{n-p+\alpha _q}\bigg)-
    \omega \bigg(\frac{2t}{n-p+1}\bigg)\bigg|}=$$
    \begin{equation}\label{27.01.11-12:38:16}
    =O(1)\frac{\alpha _q-1}{n-p+1}\omega \bigg(
    \frac{1}{n-p+1}\bigg)=\frac{O(1)}{(n-p+1)(1-q)}\omega \bigg(
    \frac{1}{n-p+1}\bigg).
    \end{equation}
Comparing (\ref{27.01.11-12:38:05}), (\ref{27.01.11-12:38:16}) and taking into account that
by (\ref{22.03.10-14:08:31}),
    \begin{equation}\label{3.03.11-20:07:08}
    {K}_{p,q}=\frac{O(1)}{(1-q)^{\delta (p)-1}},
    \end{equation}
we arrive at (\ref{25.01.11-15:00:49}).

Combining (\ref{23.03.10-19:08:49}) and (\ref{25.01.11-15:00:49}), we obtain the estimate
    $$\mathcal{E}(C^q_\beta H_\omega ;V_{n,p})\geqslant \frac{q^{n-p+1}}{\pi^2 p}
    {K}_{p,q}\int_{0}^{\pi /2}\omega \bigg(
    \frac{2t}{n-p+1}\bigg)\sin t\,dt
    +$$
    \begin{equation}\label{27.03.10-15:21:49}
    +O(1)\frac{q^{n-p+1}\omega (\pi )}{(1-q)^{\delta (p)}p(n-p+1)},
    \end{equation}
in which $O(1)$ is a quantity uniformly bounded in $n,$ $p,$ $q$, $\omega $ and $\beta $.
From (\ref{24.03.10-16:07:09}) and (\ref{27.03.10-15:21:49}) we get asymptotic formula
(\ref{27.03.10-15:23:27}). Theorem \hyperlink{15.04.11-16:12:48}{1} is
proved.\hspace{\stretch{1}}$\blacksquare$

\emph{Proof of Theorem \hyperlink{15.04.11-16:24:36}{3}.} Since the sequence $e_k(\omega )$
is monotonically decreasing (see (\ref{22.12.10-23:18:36})), from (\ref{6.05.10-14:55:11})
and (\ref{11.05.10-12:13:35}) it follows that
\begin{equation}\label{3.03.11-19:30:32}
    \mathcal{E}(C^q_\beta H_\omega ;V_{n,p})\leqslant \frac{1}{p}\sum\limits_{
    k=n-p+1}^{n}\mathcal{E}(C^q_\beta H_\omega ;S_{k-1})\leqslant $$
    $$\leqslant \frac{q^{n-p+1}}{p}\Bigg(\frac{4}{
    \pi^2 }\frac{1-q^p}{1-q}\textbf{K}(q)e_{n-p+1}(\omega )+
    \frac{O(1)}{(1-q)^{\delta (p)}(n-p+1)}\omega \bigg(\frac{1}{n-p+1}\bigg)\Bigg).
\end{equation}
On the other hand, if we analyze the proof of (\ref{23.03.10-19:08:49}), it is easy to see
that for any function $\varphi \in H_\omega $ the estimate
\begin{equation}\label{3.03.11-19:44:22}
    \mathcal{E}(C^q_\beta H_\omega ;V_{n,p})\geqslant \frac{q^{n-p+1}}{\pi p} \bigg(
    {\big|I_{n,p,q,\beta }(\varphi ) \big|+\frac{O(1)}{(1-q)^{\delta
    (p)}(n-p+1)}\|\Delta (\varphi ,\cdot) \|_C\bigg)}
\end{equation}
holds. For the function $\varphi ^*(t)$ defined by (\ref{27.01.11-11:33:54}), we have from
(\ref{27.01.11-12:38:05})--(\ref{3.03.11-20:07:08}) that
    $$I_{n,p,q,\beta }(\varphi^* )=\frac{{K}_{p,q}}{\pi }\int_{0}^{\pi /2}\omega \bigg(
    \frac{2t}{n-p+1}\bigg)\sin t\,dt+$$
    \begin{equation}\label{3.03.11-20:13:36}
    +\frac{O(1)}{(1-q)^{\delta (p)}(n-p+1)}\omega \bigg(
    \frac{1}{n-p+1}\bigg).
    \end{equation}
Since $\|\Delta(\varphi ^*,\cdot)\|_C=\frac{1}{2}\omega (\frac{\pi }{n-p+\alpha _q})$,
comparing (\ref{3.03.11-19:44:22}) and (\ref{3.03.11-20:13:36}), we obtain
\begin{equation}\label{3.03.11-20:23:12}
    \mathcal{E}(C^q_\beta H_\omega ;V_{n,p})\geqslant
    \frac{q^{n-p+1}}{p}\Bigg(\frac{K_{p,q}}{\pi ^2}e_{n-p+1}(\omega )+
    \frac{O(1)}{(1-q)^{\delta (p)}(n-p+1)}\omega \bigg(
    \frac{1}{n-p+1}\bigg)\Bigg).
\end{equation}
From (\ref{3.03.11-19:30:32}), (\ref{3.03.11-20:23:12}) and (\ref{15.02.11-11:48:33}),
relation (\ref{3.03.11-16:11:11}) follows. Theorem \hyperlink{15.04.11-16:24:36}{3} is
proved. \hspace{\stretch{1}}$\blacksquare$

\emph{Contact information}: \href{http://www.imath.kiev.ua/~funct}{Department of the Theory
of Functions}, Institute of Mathematics of Ukrainian National Academy of Sciences, 3,
Tereshenkivska st., 01601, Kyiv, Ukraine \vskip 0.2 cm \emph{E-mail}:
\href{mailto:serdyuk@imath.kiev.ua}{serdyuk@imath.kiev.ua},
\href{mailto:ievgen.ovsii@gmail.com}{ievgen.ovsii@gmail.com}

\end{document}